\documentclass{agtart_a}
\pdfoutput=1

\usepackage{amscd}
\usepackage{pinlabel}

\title{K\"ahler decomposition of 4--manifolds}

\author{R Inan\c{c} Baykur}
\givenname{R Inan\c{c}}
\surname{Baykur}
\address{Department of Mathematics\\
Michigan State University\\\newline
East Lansing MI 48824\\
USA}
\email{baykur@msu.edu}
\urladdr{}

\volumenumber{6}
\issuenumber{}
\publicationyear{2006}
\papernumber{43}
\startpage{1239}
\endpage{1265}

\doi{}
\MR{}
\Zbl{}


\received{13 May 2006}
\revised{}
\accepted{26 June 2006}
\published{}
\publishedonline{11 September 2006}
\proposed{}
\seconded{}
\corresponding{}
\editor{}
\version{}

\arxivreference{math.GT/0601396}

\makeatletter
\def\cnewtheorem#1[#2]#3{\newtheorem{#1}{#3}[section]
\expandafter\let\csname c@#1\endcsname\c@theorem}
\makeatother

\let\xysavmatrix\xymatrix
\def\xymatrix{\disablesubscriptcorrection\xysavmatrix}
\AtBeginDocument{\let\bar\wbar\let\tilde\wtilde\let\hat\what}

\newcommand{\CP}{\mathbb{CP}}

\newtheorem{theorem}{Theorem}[section]
\cnewtheorem{lemma}[theorem]{Lemma}
\cnewtheorem{proposition}[theorem]{Proposition}
\cnewtheorem{corollary}[theorem]{Corollary}

\theoremstyle{remark}
\cnewtheorem{definition}[theorem]{Definition}
\cnewtheorem{remark}[theorem]{Remark}
\cnewtheorem{example}[theorem]{Example}
\cnewtheorem{question}[theorem]{Question}
\def\dfn{\emph}

\begin{document}

\begin{asciiabstract}
In this article we show that every closed oriented smooth 4--manifold can
be decomposed into two codimension zero submanifolds (one with reversed
orientation) so that both pieces are exact Kahler manifolds with
strictly pseudoconvex boundaries and that induced contact structures on
the common boundary are isotopic. Meanwhile, matching pairs of Lefschetz
fibrations with bounded fibers are offered as the geometric counterpart
of these structures. We also provide a simple topological proof of the
existence of folded symplectic forms on 4--manifolds.
\end{asciiabstract}

\begin{abstract}
In this article we show that every closed oriented smooth 4--manifold can
be decomposed into two codimension zero submanifolds (one with reversed
orientation) so that both pieces are exact K\"ahler manifolds with
strictly pseudoconvex boundaries and that induced contact structures on
the common boundary are isotopic. Meanwhile, matching pairs of Lefschetz
fibrations with bounded fibers are offered as the geometric counterpart
of these structures. We also provide a simple topological proof of the
existence of folded symplectic forms on 4--manifolds.
\end{abstract}

\maketitle

\section{Introduction}

One possible strategy for understanding oriented smooth $4$--manifolds is
to break them up into more tractable classes of manifolds in a controlled
manner. Situated in the intersection of complex, symplectic and Riemannian
geometries, K\"ahler manifolds are the best known candidates to be pieces
of such a decomposition. The main goal of this article is to show that
this can be achieved for any closed oriented smooth $4$--manifold $X$. We
decompose $X$ into two exact K\"ahler manifolds with strictly pseudoconvex
boundaries, up to orientation, such that contact structures on the common
boundary induced by the maximal complex distributions are isotopic.

The decomposition gives rise to a globally defined $2$--form on $X$,
which we will call a \emph{(nicely) folded K\"ahler structure}, and
it belongs to a larger family of $2$--forms: \emph{folded
symplectic structures}. Cannas da Silva showed in \cite{C} that any
closed smooth oriented $4$--manifold can be equipped with a folded
symplectic form, by using a version of the h-principle defined for
folding maps by Eliashberg. In \fullref{SimpleFoldedSymplectic}, we
introduce a way to construct some simple examples of folded symplectic
$4$--manifolds. Afterwards we reprove the existence fact by constructing a
folded symplectic form $\omega$ for a given handlebody decomposition of
$X$, essentially by means of simple Kirby calculus and contact topology
(\fullref{AllFolded}). The main ingredient there will be achiral
Lefschetz fibrations, and recent work of Etnyre and Fuller \cite{EF}
will play a key role in our construction.

Next, we switch gears, and using several results on compact Stein
surfaces and Lefschetz fibrations with bounded fibers (mainly Harer \cite{Ha},
Eliashberg \cite{El1}, Gompf \cite{Go2}, Loi--Piergallini \cite{LP},
and Akbulut--Ozbagci \cite{AO}), we prove the aforementioned
decomposition theorem (\fullref{SteinDecom}). In fact we obtain
a stronger result, as the pieces of this decomposition are actually
Stein manifolds with strictly pseudoconvex boundaries. It was first
shown by Akbulut and Matveyev in \cite{AM} that a closed oriented smooth
$4$--manifold $X$ could always be decomposed into Stein pieces, but there
was no particular information one could use to argue for matching the
induced contact structures on the separating hypersurface. Our proof
follows an alternative way via open book decompositions, and we conclude
that the Stein structures can be chosen to agree on the common contact
boundary.

In the last section, we introduce folded K\"ahler structures, and discuss
some properties they enjoy, after showing that all closed oriented
smooth $4$--manifolds admit them (\fullref{AllKahler}). This improves
the folded symplectic existence result, and indeed both structures
we construct are shown to be equivalent on the symplectic level. The
collection of these discussions leads us to define \emph{folded Lefschetz
fibrations}, roughly speaking, pairs of positive and negative Lefschetz
fibrations over disks with bounded fibers which agree on the common
boundary through induced open book decompositions. We prove that any
nicely folded K\"ahler $4$--manifold, possibly after an orientation
preserving diffeomorphism, admits compatible folded Lefschetz fibrations
(\fullref{CompatibleFoldedLF}).

\subsubsection*{Acknowledgments}
The author would like to thank Ron Fintushel for
his support and for helpful conversations. The author would also like
to thank Burak Ozbagci for commenting on the draft of this paper. The
author was partially supported by NSF Grant DMS0305818.

\section{Preliminaries}

\subsection{Lefschetz fibrations}
All manifolds and maps in this article are assumed to be smooth. A
\dfn{Lefschetz fibration} on an oriented 4--manifold $X$, possibly
with boundary, is a map $f\colon\, X\to \Sigma$, where $\Sigma$ is a
compact oriented surface, such that each critical point of $f$ lies
in the interior of $X$ and has a local model $f(z_1,z_2)=z_1 z_2$,
given by orientation preserving charts both on $X$ and $\Sigma$. These
singularities are obtained by attaching $2$--handles to regular fibers
with framing $-1$ with respect to the framing induced by the fiber. We
will refer to these $2$--handles as \dfn{positive Lefschetz handles}. An
\dfn{achiral Lefschetz fibration} is defined the same way, except that the
given charts around critical points are allowed to reverse orientation. In
other words, the $2$--handles can be glued with framing $+1$ with respect
to fiber framing, too. Also recall that a Lefschetz pencil is a map $f\co
X \setminus \{b_1,\ldots,b_k\}\to S^2$, such that around any \emph{base
point} $b_i$ it has a local model $f(z_1,z_2)=z_1/z_2 \,$, preserving
the orientations, and that $f$ is a Lefschetz fibration elsewhere. An
\dfn{achiral Lefschetz pencil} is then defined by allowing orientation
reversing charts around the base points as well. Critical points or base
points with orientation reversing charts are called \dfn{negative critical
points} or \dfn{negative base points}, whereas the other critical points
or base points are \dfn{positive}. For a detailed treatment of this
topic and proofs of some facts quoted below, the reader is advised to
turn to Gompf and Stipsicz \cite{GS}.

A Lefschetz fibration is said to be \dfn{allowable} if all its vanishing
cycles are homologically nontrivial in the fiber. Particularly, we will
be interested in allowable Lefschetz fibrations over $D^2$ with bounded
fibers. In the literature, this type of Lefschetz fibration having only
positive critical points is called a \dfn{PALF}. Similarly, when the
critical points are instead all negative, we will call the fibration
a \dfn{NALF}.

Given a compact oriented genus $g$ surface $F$ with $m$ boundary
components and $r$ marked points on it, the \dfn{mapping class group\,}
$\Gamma_{g,m}^r$ is defined as the group of orientation preserving
self-diffeomorphisms of $F$ fixing marked points and $\partial F$
pointwise, modulo isotopies of $F$ fixing marked points and $\partial F$
pointwise. It can be shown that $\Gamma_{g,m}^r$ is generated by positive
(right-handed) and negative (left-handed) Dehn twists. Importantly,
isotopy type of a surface bundle over $S^1$ with fiber closed oriented
surface $F$ is determined by the return map of a flow transverse to
the fibers, which can be identified with an element $\mu \in \Gamma_g$,
called \dfn{monodromy} of this fibration.

Let $f\co X \to D^2$ be an achiral Lefschetz fibration, where the
regular fiber $F$ is an oriented genus $g$ surface with $m$ boundary
components, and suppose all critical points of the fibration lie on
various fibers. Select a regular value $p$ in the interior of $D^2$,
an identification of $f^{-1}(p)\cong F$, and a collection of arcs $s_1,
\ldots, s_k$ in the interior of $D^2$ with each $s_i$ connecting $p$
to a distinct critical value, and all disjoint except at $p$. We index
the critical values as well, so that each arc $s_i$ is connected to a
critical value $q_i$ and that they appear in a counterclockwise order
around the point $p$. Now if we take a regular neighborhood of each arc
away from remaining critical points and consider the union of these,
we obtain a disk $V$ and an $F$--bundle over $\partial V = S^1$. The
monodromy of this fibration is an element $\mu \in \Gamma_{g,m}$, which
is called the \dfn{global monodromy} of the achiral Lefschetz fibration
$f$. It is well-known that this data gives a handlebody description of
$X$, and vice versa. We call the ordered set of arcs $\{s_1, \ldots,
s_k\}$ a \dfn{representation} of the achiral Lefschetz fibration $f$.

Next is a standard fact which was first observed by Harer:

\begin{theorem}[Harer \cite{Ha}]
\label{Harer}
Let $X$ be a $4$--manifold with boundary. Then $X$ admits an achiral
Lefschetz fibration over $D^2$ with bounded fibers if and only if
it admits a handlebody decomposition with no handle of index greater
than two.
\end{theorem}


\subsection{ Open book decompositions}

An \dfn{open book decomposition} of a $3$--manifold $M$ is a pair $(B,
f)$ where $B$ is an oriented link in $M$, called the \dfn{binding},
and $f\co M\setminus B \to S^1$ is a fibration such that  $f^{-1}(t)$ is
the interior of a compact oriented surface $F_t \subset M$ and $\partial
F_t=B$ for all $t \in S^1$. The surface $F=F_t$, for any $t$, is called
the \dfn{page} of the open book. The \dfn{monodromy} of an open book is
given by the return map of a flow transverse to the pages and meridional
near the binding, which is an element $\mu \in \Gamma_{g, m}$, where $g$
is the genus of the page $F$, and $m$ is the number of components of $B =
\partial F$.

Suppose we have an achiral Lefschetz fibration $f\co X \to D^2$ with
bounded regular fiber $F$, and let $p$ be a regular value in the interior
of the base $D^2$. Composing $f$ with the radial projection $D^2 \setminus
\{p\} \to \partial D^2$ we obtain an open book decomposition on $\partial
X$ with binding $\partial f^{-1}(p)$. Identifying $f^{-1}(p) \cong F$,
we can write $\partial X=(\partial F\times D^2)\cup f^{-1}(\partial
D^2)$. Thus we view $\partial F\times D^2$ as the tubular neighborhood of
the binding $B=\partial f^{-1}(p)$, and the fibers over $\partial D^2$ as
its \dfn{truncated pages}. The monodromy of this open book is prescribed
by that of the achiral fibration \cite{Ha}. In this case, we say the open
book $(B, f|_{\partial X \setminus B})$ \dfn{bounds} or \dfn{is induced
by} the achiral Lefschetz fibration $f\co X\to D^2$. Recalling that any
closed oriented $3$--manifold can be bounded by a $4$--manifold with only
$0$--, $1$-- and $2$-- handles, it is fairly easy to see that any open book
decomposition bounds such an achiral Lefschetz fibration over a disk.

We would like to describe an elementary modification of these structures:
Let $f\co X \to D^2$ be an achiral Lefschetz fibration with bounded
regular fiber $F$. Attach a $1$--handle to $\partial F$ to obtain $F'$,
and then attach a positive (resp. negative) Lefschetz $2$--handle along
an embedded loop in $F'$ that goes over the new $1$--handle exactly
once. This is called a \dfn{positive stabilization} (resp.\,\dfn{negative
stabilization}) of $f$. A positive (resp. negative) Lefschetz handle is
attached with framing $-1$ (resp.\,$+1$) with respect to the fiber,
and thus it introduces a positive (resp. negative) Dehn twist on
$F'$. If the focus is on the $3$--manifold, one can totally forget
the bounding $4$--manifold and view all the handle attachments in the
$3$--manifold. Either way, stabilizations correspond to adding canceling
handle pairs, so diffeomorphism types of the underlying manifolds do
not change, whereas the achiral Lefschetz fibration and the open book
decomposition change in the obvious way. It turns out that stabilizations
preserve more than the underlying topology, as we will discuss shortly.


\subsection{ Contact structures and compatibility }

A $1$--form $\alpha \in \Omega^1(M)$ on a $(2n{-}1)$--dimensional oriented
manifold $M$ is called a \dfn{contact form} if it satisfies $\alpha \wedge
(d\alpha)^{n-1} \neq 0$. An \dfn{oriented contact structure} on $M$ is
then a hyperplane field $\xi$ which can be globally written as kernel
of a contact $1$--form $\alpha$. In dimension three, this is equivalent
to asking that $d\alpha$ be nondegenerate on the plane field $\xi$.

A contact structure $\xi$ on a $3$--manifold $M$ is said to be
\dfn{supported by an open book} $(B,f)$ if $\xi$ is isotopic to a
contact structure given by a $1$--form $\alpha$ satisfying $\alpha>0$ on
positively oriented tangents to $B$ and $d\alpha$ is a positive volume
form on every page. When this holds, we say that the open book $(B,f)$
is \dfn{compatible with the contact structure} $\xi$ on $M$.

Improving results of Thurston and Winkelnkemper \cite{TW}, Giroux proved
the following groundbreaking theorem regarding compatibilty of open
books and contact structures:

\begin{theorem} [Giroux \cite{Gi2}] \label{Giroux}
Let $M$ be a closed oriented $3$--manifold. Then there is a one-to-one
correspondence between oriented
contact structures on $M$ up to isotopy and open book decompositions of
$M$ up to positive stabilizations and isotopy.
\end{theorem}

Considering contact $3$--manifolds as boundaries of certain $4$--manifolds
together with some compatibility conditions is a current focus of research
in low dimensional topology. From the contact topology point of view,
it is the study of different types of \dfn{fillings} of a fixed contact
manifold. In dimension four, there are essentially two considerations,
yet we formulate them for all dimensions: Let $(X^{2n},\omega)$ be
a symplectic manifold with cooriented nonempty boundary $M=\partial
X$. If there exists a \emph{Liouville vector field} (aka \emph{symplectic
dilation}) $\nu$ defined on a neighborhood of $\partial X$ pointing out
along $\partial X$, then we obtain a positive contact structure $\xi$
on $\partial X$, which can be written as the kernel of contact $1$--form
$\alpha = \iota_{\nu} \omega|_{\partial X}$. When this holds, we say
$(M, \xi)$ is the \dfn{$\omega$--convex boundary} or \dfn{strongly convex
boundary} of $(X, \omega)$. For the sake of entirety, note when $\nu$
points inside, we obtain a negative contact structure instead, and in
this case we say $(M, \xi)$ is the \dfn{$\omega$--concave boundary of
$(X, \omega)$}.

Now if $(X^{2n}, J)$ is almost-complex, then the complex tangencies on
$M=\partial X$ give a uniquely defined oriented hyperplane field. It
follows that there is a $1$--form $\alpha$ on $M$ such
that $\xi = Ker \alpha$. We define the \dfn{Levi form} on $M$ as
$d\alpha|_{\xi}(\cdot, J \cdot)$. If this form is positive definite then
$(M, \xi)$ is said to be \dfn{strictly $J$--convex boundary of $(X,J)$},
and if it is $J$--convex for an unspecified $J$ (for instance when $J$ is
tamed by a given symplectic form), we say $(M, \xi)$ is \dfn{strictly
pseudoconvex boundary}. If $(X, \omega, J)$ is an almost-K\"ahler
manifold, i.e. a manifold equipped with a symplectic form $\omega$ and
a compatible almost-complex structure $J$, then it can be shown that
strict pseudoconvexity of the boundary is equivalent to the condition
that ${\omega|}_{\xi } > 0$ in dimenson $2n=4$. We would like to remark
that all these definitions can be formulated in more generality for
hypersurfaces in $X^{2n}$, not necessarily for $\partial X$ only.

For detailed and comparative discussions of these concepts, as well
as proofs of some facts mentioned in the next subsection, the reader
can turn to Eliashberg--Gromov \cite{ElGr} and Etnyre \cite{E}. Also
for further basic notions from contact topology of $3$--manifolds such
as Legendrian knots, Thurston--Bennequin framing, or convex surfaces,
which we will occasionally use in this paper, see for example
Ozbagci--Stipsicz \cite{OS}.


\subsection{ Stein manifolds}

A smooth function $\psi\co X \to \R$ on a complex manifold $X$ of real
dimension $2n$ is called \dfn{strictly plurisubharmonic} if $\psi$
is strictly subharmonic on every holomorphic curve in $X$. We call
a complex manifold $X$ \dfn{Stein}, if it admits a proper strictly
plurisubharmonic function $\psi\co X \to [0, \infty)$ (\,after Grauert
\cite{Gra}). Thus a compact manifold $X$ with boundary which is equipped
with a complex structure in its interior is called \dfn{compact Stein}
if it admits a proper strictly plurisubharmonic function which is constant
on the boundary.

Given a function $\psi\co X \to \R$ on a Stein manifold, we can define
a $2$--form $\omega_{\psi} = -dJ^* d\psi$. It turns out that $\psi$
is a strictly plurisubharmonic function if and only if the symmetric
form $g_{\psi}(\cdot, \cdot) = \omega_{\psi}(\cdot, J \cdot)$ is
positive definite. So every Stein manifold $X$ admits a K\"ahler
structure $\omega_{\psi}$, for any strictly plurisubharmonic function
$\psi\co X \to [0, \infty)$. It is easy to see that the restriction of
$\omega_{\psi}$ to each level set $\psi^{-1}(t)$ gives a Levi form on
$\psi^{-1}(t)$, implying that all nonsingular level sets of $\psi$ are
strictly pseudoconvex hypersurfaces. Thus in this article, we equivalently
call a Stein manifold a \dfn{strictly pseudoconvex manifold}. Moreover,
it was observed in \cite{ElGr} that the gradient vector field of $\psi$
defines a (global) Liouville vector field $\nu=\nabla_{\psi}$, making all
nonsingular level sets $\omega_{\psi}$--convex. Hence, Stein manifolds
exhibit strongest filling properties for a contact manifold which can
be realized as their boundary.

In this article, we are mainly interested in compact Stein
surfaces. Another characterization of these manifolds, which might be
called ``the topologist's fundamental theorem of compact Stein surfaces'',
is due to Eliashberg, and was made explicit by Gompf in dimension four:

\begin{theorem} [Eliashberg \cite{El1}, Gompf \cite{Go2}]
\label{Eliashberg}
A smooth oriented compact $4$--manifold with boundary is a
Stein surface, up to orientation preserving diffeomorphisms, if and
only if it has a handle decomposition $X_0 \cup h_1 \cup \ldots \cup h_m
$, where  $X_0$ consists of \,$0$-- and \,$1$--handles and each $h_i$,
$1\leq i \leq m$, is a $2$--handle attached to $X_i= X_0 \cup h_1 \cup
\ldots \cup h_i$ along a Legendrian circle $L_i$ with framing $tb(L_i) -1 $.
\end{theorem}

All structures we have introduced so far meet in the following theorem:

\begin{theorem} [Loi--Piergallini \cite{LP}, also see Akbulut--Ozbagci \cite{AO}]
\label{PALF}
An oriented compact $4$--manifold with boundary is a Stein
surface, up to orientation preserving diffeomorphisms, if and only if
it admits a PALF.
\end{theorem}

Throughout the article, we give ourselves the freedom of using the
prefix `anti' as a shorthand, whenever an oriented manifold $X$ admits a
structure when the orientation on $X$ is reversed; like anti-symplectic,
anti-K\"ahler, or anti-Stein. For Lefschetz fibrations and open books
though, we use `positive' and `negative' adjectives to distinguish two
possible cases.


\section{Simple folded symplectic structures}
\label{SimpleFoldedSymplectic}

The definition of symplectic (or anti-symplectic) structures can be
enlarged as follows in order to cover a larger family of manifolds,
which was shown by Cannas da Silva \cite{C} to contain entire family of
closed oriented smooth $4$--manifolds:

\begin{definition} \label{FoldedSymplectic}
A \dfn{folded symplectic form} on a smooth $2n$--dimensional manifold $X$
is a closed $2$--form $\omega$ such that $\omega^n$ is transverse to the
$0$--section of $\Lambda^{2n} T^*X$, and whenever this intersection
is nonempty, $\omega^{n-1}$ does not vanish on the hypersurface
$H=(\omega^n)^{-1}(0)$, called the \dfn{fold}.
\end{definition}

For an oriented $X$, the kernel of $\omega$ on $H$ integrates to
a foliation called \dfn{null-foliation}. Martinet's singular form
$x_1 dx_1 \wedge dy_1 + dx_2 \wedge dy_2 + \cdots + dx_n \wedge dy_n$ on
$\R^{2n}$ defines the standard folded symplectic structure, as every
folded symplectic form can be expressed in this way in an appropriate
Darboux coordinate system around any point on the fold. There is also
a simple folded structure that every even dimensional sphere carries:
We think of $S^{2n}$ sitting in $\R^{2n+1}$, then pull back the standard
symplectic form $dx_1 \wedge dy_1 + \cdots + dx_n \wedge dy_n$ on the unit
disk bounded by the equator in $\R^{2n}$ to $S^{2n}$ by the projection
maps along the last coordinate, and finally glue them along the fold
$S^{2n-1}$ to obtain $\omega_0$. This is equivalent to doubling the
unit disk equipped with its standard symplectic form (by reversing the
orientation on one of the disks). We call this form the \dfn{standard
folded symplectic form on $S^{2n}$}.

For more on folded symplectic structures, the reader is referred to the
work of Cannas da Silva, Guillemin and Woodward \cite{C,CGW}. Here we
only consider these forms on Riemann surfaces and compact 4--manifolds,
possibly with boundaries. For the former class, folded symplectic
forms form an open and dense set in the space of $2$--forms, whereas
in dimension four openness remains but the nonvanishing condition
implies that they are nongeneric. We say an embedded surface $\Sigma
\subset X^4$ is a \dfn{folded symplectic submanifold} of $(X, \omega)$
if $\omega|_{\Sigma}$ is a folded symplectic form on $\Sigma$. Observe
that $S^2$ equipped with the standard form obtained by pulling back
$dx_1 \wedge dy_1$ embeds as a folded symplectic submanifold of $S^4$
with the standard folded symplectic form defined as the pullback of
$dx_1 \wedge dy_1 + dx_2 \wedge dy_2$ as above.

The following proposition provides several examples of folded symplectic
$4$--manifolds:

\begin{proposition} \label{SimpleModel}
Let $X$ be a closed oriented smooth $4$--manifold and $\Sigma$ be a
closed oriented surface. If $f\co X \to \Sigma^2$ is an achiral Lefschetz
fibration such that the regular fiber is a closed oriented surface $F$
which is nonzero in $H_2(X;\R)$, then $X$ admits a folded symplectic
structure $\omega$ such that fibers are symplectic and the fold $H$ is
an $F$--bundle over $S^1$. The fold $H$ splits $X$ into pieces $X_+$ and
$X_-$, and $f$ induces symplectic Lefschetz fibrations on $(X_+, \omega
|_{X_+})$ and on $(-X_-, \omega |_{X_-})$, respectively. Furthermore, any
finite set of sections can be made folded symplectic for an appropriate
choice of $\omega$. This form is canonical up to deformation equivalence
of folded symplectic forms.
\end{proposition}

We will call this type of folded symplectic structures \dfn{simple}
(after Thurston \cite{T}). Base spaces of the fibrations defined on $X_+$
and $-X_-$ are determined by an arbitrary splitting $\Sigma = \Sigma_+
\cup \Sigma_-$. Here we take $\Sigma_-=D^2$ for simplicity. Observe
that the fibration induces an exact sequence $$\pi_1(F)\to\pi_1(X)\to
\pi_1(\Sigma)\to \pi_0(F)\to 0$$ It follows that fibers are connected
if the base is simply-connected. Otherwise we can define a new achiral
Lefschetz fibration from $X$ to the finite cover of $\Sigma$ corresponding
to the finite-index subgroup $f_{\#} (\pi_1(X))$ in $\pi_1(\Sigma)$,
which has connected fibers. Finally, one can perturb $f$ to get a
fibration which has at most one critical point on each fiber. Hence,
without loss of generality, we will assume that the fibers of $f$ are
connected and critical values are distinct.


\begin{proof} [Proof of \fullref{SimpleModel}]
Start by connecting all negative critical points in the base by an
embedded arc in the complement of positive critical points, and cover it
by the images of orientation reversing charts so that we get a closed
neighborhood $\Sigma_-\cong D^2$ of this arc away from the positive
critical points. This can be done because around the regular points we
have freedom to take charts of either orientation. After we reverse the
orientation on $f^{-1}(\Sigma_-)$, the map $f\co f^{-1}(\Sigma_-) \to
\Sigma_-$ defines a negative Lefschetz fibration. Set $\Sigma_+=\Sigma
\setminus \Sigma_-$, $C=\Sigma_+\cap \Sigma_-$, $X_+=f^{-1}(\Sigma_+)$,
$X_-=f^{-1}(\Sigma_-)$, and $H=f^{-1}(C)$. If there are no negative
critical points, we can choose $\Sigma_-$ as a small disk around a
regular value which does not contain any critical values. Now let $\beta$
be a folded symplectic form on $\Sigma$ which folds over $C$, such that
it is a positive area form on $\Sigma_+$ and a negative area form on
$\Sigma_-$. These forms always exist: For example take $S^2$ with its
standard folded form $\omega_0$, and suppose $\Sigma_{\pm}$ has genus
$g_{\pm}$. Symplectic connect sum the upper-hemisphere of $S^2$ with a
closed genus $g_+$ surface equipped with a positive symplectic form,
and the lower-hemisphere with a closed genus $g_-$ surface equipped
with a negative symplectic form. This yields a folded symplectic form
on $\Sigma$, folded along $C$.

We will construct a folded symplectic form on $X$ by mimicking
Gompf's proof which generalizes Thurston's result for symplectic
fibrations to symplectic Lefschetz fibrations (see Thurston \cite{T}
and Gompf--Stipsicz \cite{GS}). Let $\zeta$ be a closed 2--form on $X$
which evaluates positively on any closed surface contained in a fiber
with the induced orientation. (We have not made any assumptions on the
type of vanishing cycles, so one might have more than one closed surface
on a fiber if there are separating vanishing cycles.) First we wish to
define a closed 2--form $\eta$ on all over $X$ which is symplectic on
each $F_y=f^{-1}(y)$, for all $y\in \Sigma$.

Let $A$ be a tubular neighborhood of $C$ in $\Sigma$ which does not
contain any critical values. Choose disjoint open balls $U_{+,k}$
around each positive and $V_{-,l}$ around each negative critical point
so that these sets do not intersect $f^{-1}(A)$ in $X$ and that in
appropriate charts the fibration map can be written as $f(z_1,z_2)=z_1
z_2$ and $f(z_1,z_2)=\wbar{z}_1 z_2$, respectively. Take the
standard forms $$\omega_{+,k}=dx_1\wedge dy_1 + dx_2 \wedge dy_2=
-\tfrac{i}{2}dz_1\wedge d\wbar{z}_1 -\tfrac{i}{2}dz_2 \wedge
d\wbar{z}_2$$ on $U_{+,k}$ and $$\omega_{-,l}=-dx_1\wedge dy_1 + dx_2
\wedge dy_2= \tfrac{i}{2}dz_1\wedge d\wbar{z}_1 -\tfrac{i}{2}dz_2 \wedge
d\wbar{z}_2$$ on $V_{-,l}$ for all $k,l$. For any $y\in f(U_{+,k})$,
$F_y\cap U_{+,k}$ is a $J_{+,k}$--holomorphic curve, where $J_{+,k}$ is
an almost-complex structure compatible with $\omega_{+,k}$. Similarly
for any $y\in f(V_{-,l})$, $F_y\cap V_{-,l}$ is $J_{-,l}$--holomorphic
curve, where $J_{-,l}$ is an almost-complex structure compatible with
$\omega_{-,l}$. Having expressed $\omega_{+,k}$ and $\omega_{-,l}$ in
terms of K$\ddot{\rm{a}}$hler forms, we can take these almost-complex
structures as $(i,i)$ and $(-i,i$), respectively. It follows that
$\omega_{+,k}|_{F_y\cap U_{+,k}}$ is symplectic, so we can extend it to a
symplectic form $\omega_y$ on the entire fiber and get $\omega_y$ defined
for all points in each $f(U_{+,k})$ this way. Do the same for all points
in $f(V_{-,l})$, for every $j$. Finally, for all remaining $y\in \Sigma$
take any symplectic form $\omega_y$ on the fiber, and rescale every
$\omega_y$ we have defined away from all $U_{+,k}$ and $V_{-,l}$ so that
they are in the same cohomology class as the restriction of $\zeta$ to
each $F_y$. Next, cover $\Sigma$ with finitely many balls $B_s$ containing
at most one critical value, and whenever they do contain a critical value,
assume they are centered at that point. Reindex $U_{+,k}$ and $V_{-,l}$,
and shrink them if necessary to make sure they lie in $f^{-1}(B_s)$
for some $s$. Define $\eta_s$ on each $f^{-1}(B_s)$ as the pullback of
$\omega_{+,s}$, $\omega_{-,s}$, or $\omega_y$ by $r_s$, where $r_s$ is the
retraction of $f^{-1}(B_s)$ to the fiber $F_y$ over the center of $B_s$,
or the union of $F_y$ either with closure of $U_{+,s}$ or with closure
of $V_{-,s}$, whenever $B_s$ contains a positive or negative critical
value, respectively. Now we can glue these forms to construct the 2--form
$\eta$ we wanted, by using a partition of unity and that each $\eta_s$
is cohomologous to $\zeta|_{f^{-1}(B_s)}$ as in Gompf--Stipsicz \cite{GS}.

We claim that $\omega_{\kappa}= \kappa \eta + f^*(\beta)$ is a folded
symplectic form on $X$, where $\kappa$ is a small enough positive
real number. $\omega_{\kappa}$ is clearly closed and symplectic in
the fiber direction. It follows that for any noncritical point $x \in
F_y$, $T_xX=T_xF_y \oplus (T_xF_y)^{\perp \eta}$. Here $f^*(\beta)$
is nondegenerate over $(T_xF_y)^{\perp \eta}$ for all $x\notin H$,
implying that for sufficiently small $\kappa$, $\omega_{\kappa}$ is
nondegenerate on $X\setminus (H \, \bigcup_s (U_{+,s} \cup V_{-,s}))$. On
the other hand $\omega_{\kappa} |_{U_{+,s}} = \kappa \omega_{+,s} +
f^* (\beta)$ and $\omega_{\kappa} |_{V_{-,s}} = \kappa \omega_{-,s}
+ f^* (\beta)$. Therefore for any nonzero $v \in TU_{+,s}$, we have
$$\omega(v,J_{+,s}\, v) = \kappa g(v, v)_{+,s}^2 + \beta(f_*(v), i f_*(v))
>0 \, ,$$ where $g(-,-)_{+,s}$ is the metric induced from $\omega_{+,s}$
and $J_{+,s}$. Likewise, for any nonzero $v \in TV_{-,s}$, we will have
$$\omega(v,J_{-,s} \, v) = \kappa g(v, v)_{-,s}^2 + \beta(f_*(v), -i
f_*(v)) >0,$$ $g(-,-)_{-,s}$ being the metric induced from $\omega_{-,s}$
and $J_{-,s}$. (Recall that $\beta$ is negative on $\Sigma_-$.) Hence
$\omega_{\kappa}$ is symplectic everywhere on $X$ except $H$, where it
vanishes transversely. Moreover, $f^*(\beta)$ is a folded symplectic form
on any section, so taking $\kappa$ even smaller, we can as well assume
that any finite collection of sections of $f$ are folded symplectic. It
is easy to check that the folded symplectic form we get satisfies all
the other declared properties. (Also see \fullref{AlternativeProof}).
\end{proof}

The homological assumption in the theorem is a very mild one. If $S$ is
the set of critical points of the achiral fibration $f\co X\to\Sigma$,
then the tangencies of the fibers define a complex line bundle $L=Ker(df)$
on $X \setminus S$, which extends uniquely over $X$. It follows that
unless we have a torus fibration, the regular fiber $F$ is essential,
since $\langle c_1(L), F\rangle=\chi(F)$. Also if the fibration is obtained from a
pencil by blowing up the base points, the exceptional spheres will become
sections of the fibration, guaranteeing that the fibers are essential
in the homology.

\begin{remark} \label{AlternativeProof} \rm
Alternatively, the folded symplectic form in \fullref{SimpleModel} can
be constructed by using the \dfn{folding} operation described by Cannas
da Silva, Guillemin and Woodward \cite{CGW}. Restrictions of $\beta$ on
$\Sigma_+$ and on $\Sigma_-$ give well-defined area forms $\beta_+$ and
$\beta_-$, respectively. Gompf's method can be used to define a symplectic
form $\kappa_+ \eta + f^*(\beta_+)$ on $X_+$, where $\eta$ is a 2--form
on $X$ that restricts to the fibers as a (positive) symplectic form and
$\kappa_+$ is a small enough positive real number. The orientation on
the base together with the orientation on the regular fiber determines
the orientation of the total space, and thus by taking the opposite
orientation on $\Sigma_-$ but keeping the orientation on $F$, one orients
$-X_-$. Let $\wbar{f}\co -X_- \to \Sigma_-$ be the fibration defined
by taking orientation-preserving charts for $f\co X_-\to \Sigma_-$, then
we can define a symplectic form $\kappa_- \eta + \wbar{f}^*(-\beta_-)$
on $-X_-$ (as $-\beta_-$ is the area form on $\Sigma_-$) by following
the same construction method. Observe that $\wbar{f}^*(-\beta_-)
= f^*(\beta_-)$. Hence, setting $\kappa = min\{\kappa_+, \kappa_-\}$,
we obtain two symplectic manifolds $(X_+, \omega_+)$ and $(-X_-,
\omega_-)$, where $\omega_{\pm}=\kappa \eta + f^*(\beta_{\pm})$. Let
$\iota_{\pm}$ be the inclusions of boundaries into $\pm X_{\pm}$,
then $\iota_+^*(\omega_+)= \kappa \eta = \iota_-^*(\omega_-)$ and the
orientations of both null-foliations agree. Thus we can glue these
pieces to obtain a folded symplectic structure on $X_+ \cup X_- =X$,
which agrees with $\omega_+$ and $\omega_-$ in the complement of a tubular
neighborhood of the fold $\partial X_+=H= -\partial X_-$. (See \cite{CGW}
for details.) This form is deformation equivalent to the form $\kappa
\eta + f^*(\beta)$ in \fullref{SimpleModel}.
\end{remark}


\section{Existence of folded symplectic structures on closed oriented
4--manifolds} \label{FoldedExistence}

Here we show that any closed oriented smooth 4--manifold $X$ can be
equipped with a folded symplectic form. For the sake of completeness,
we start by outlining Etnyre and Fuller's proof \cite{EF} that every
4--manifold admits an achiral Lefschetz fibration after a surgery along
a framed circle: Take a handlebody decomposition of $X$ with one
$0$-- and one $4$--handle, let $X_1$ denote the union of the $0$--handle,
$1$--handles and $2$--handles, and $X_2$ denote the union of the $3$--handles
and the $4$--handle. By \fullref{Harer} there exist achiral Lefschetz
fibrations $f_i\co X_i\to D^2$, which necessarily have bounded fibers,
and stabilizing both fibrations we may as well assume the fibers have
connected boundaries. After a possible slight modification of the
handlebody decomposition, Etnyre and Fuller manipulate the contact
structures on the boundaries so that they are both overtwisted and
homotopic as plane fields. Then it follows from results of Eliashberg
and Giroux that we have isotopic contact structures, and thus the
induced open books are the same, possibly after some stabilizations and
isotopies. Denoting the final manifolds and fibrations with $X_i$ and
$f_i$ again, we may therefore assume that the open book decompositions
induced by these fibrations on the common boundary $H=\partial
X_1=-\partial X_2$ are the same, so we can glue both pieces of $X$
back along the truncated pages, and obtain an achiral Lefschetz fibration
\[
f_1\cup f_2\co W= X_1\bigcup_{f_1^{-1}(\partial D^2)=f_2^{-1}(\partial
D^2)} X_2 \longrightarrow S^2.
\]
To recover $X$ we need to glue $S^1\times D^2_1$ to $S^1\times D^2_2$,
where $$S^1\times D_i^2=\partial X_i\setminus f_i^{-1}(\partial D^2).$$
Filling the boundary of $W$ with an $S^1\times D^3$ gives the same result,
so we can view $W$ as $X\setminus N$ where $N$ is a neighborhood of an
embedded curve $\gamma \subset X$. Now, if we instead add on a $D^2\times
S^2$ so that each $\partial D^2 \times \{pt\}$ is identified with $S^1
\times \{pt\}$, we can extend the fibration on $W$ by the projection
on the $S^2$ component of $D^2\times S^2$. Hence, we obtain an achiral
Lefschetz fibration over $S^2$ on the resulting manifold $Y$, where the
section $S$ of this fibration discussed by Etnyre and Fuller \cite{EF}
can be taken as $0\times S^2$ coming from the glued in $D^2\times S^2$,
implying $S$ has trivial normal bundle in $Y$.

We will refer the following as the \dfn{standard model}: Consider
$S^4$ with the standard folded symplectic structure $\omega _0$
described before, and take $S^4\cap\{x_4=0\}$ vertical to the fold $H_0
= S^4\cap\{x_5=0\}$. Take \mbox{$S_0=S^4\cap\{x_4=0=x_3\}\cong S^2$}
which intersects the fold along the circle $C_0=\smash{\bigl\{x_1^2+x_2^2=1
~|~x_3=x_4=x_5=0\bigr\}}$. It is easy to see that $\omega _0$ restricts to this
$S_0$ as the standard folded symplectic form on $S^2$, folded along $C_0$,
and symplectic on the normal disks to $S_0$. Fix a disk neighborhood $M_0$
of $S_0$ so that $\omega_0$ evaluates as $1$ on each normal disk. That
is, each normal disk projects onto unit disk $\bigl\{x_3^2+x_4^2\leq1
\,|x_1=x_2=x_5=0\bigr\}$ symplectomorphically. By restricting  $\omega_0$,
we get two folded symplectic manifolds $M_0 \cong S^2 \times D^2$ and
$N_0 = S^4 \setminus M_0 \cong D^3 \times S^1$, with folds $S^1 \times
D^2$ and $D^2 \times S^1$, respectively.

The existence of the section $s\co S^2\to S\subset X$ guarantees that the
fiber of the achiral Lefschetz fibration $f\co Y \to S^2$ is homologically
essential and therefore there exists a folded symplectic form $\omega$
as described in \fullref{SimpleModel}. This restricts to
$Y\setminus M$, where $M\cong S^2\times D^2$ is a neighborhood of $S$. We
may assume $\omega$ is constructed such that $M$ is identified with $M_0$
in the standard model above as follows: Let $\phi\co M\to M_0$ be an
orientation preserving diffeomorphism such that $\phi$ is orientation
preserving on the spheres (and on the normal directions as well), and
that it maps the upper-hemisphere of $S_0$ (where $\omega_0$ is positive)
to the positive part of $S$. Then one can start the construction in the
proof of \fullref{SimpleModel} with the folded symplectic form
$s^*  \phi^* (\omega_0)$ on the base sphere, which naturally restricts
to an area form on each hemisphere. We can also modify the symplectic
form $\kappa \eta$ on the fibers so that it is symplectomorphic to
$\phi^*(\omega_0)$ on the normal disks to $S$, each of which lies on
a fiber.

Hence we obtain a folded symplectic form $\omega$ on
$X$ such that $(M,\omega|_M)$ is folded symplectomorphic to
$(M_0,\omega_0|_{M_0})$. This allows us to trade $M$ for $N\cong S^1\times
D^3$ and extend the folded symplectic structure to $(Y\setminus M)
\cup N \cong X$. The effect of this surgery on the fold of $Y$ is to
turn the surface fibration over $S^1$ into an open book decomposition
on the resulting fold. The core curve of $N$ sits in the $3$--manifold
as the binding of this open book and therefore it carries a canonical
framing. We have proved:

\begin{theorem} \label{AllFolded}
Every closed oriented smooth $4$--manifold $X$ admits a folded symplectic
structure. Furthermore, there exist folded symplectic forms on $X$ with
connected folds, such that a surgery along a framed curve which lies in
the fold results in a simple folded symplectic manifold.
\end{theorem}

\begin{remark} \rm
Away from the framed curve $\gamma$ in $X$, the folded symplectic model
we have constructed is the restriction of the simple model discussed
in the previous section, and as we will see shortly, the pieces are
Stein and anti-Stein. So for any sort of pseudo-holomorphic curve
counting with respect to this folded symplectic structure, the focus
would be understanding the limit behaviors around $\gamma$ of the
curves in the moduli space, where we do have a standard model, namely
$(N_0 , \omega_0|_{N_0})$ above. (For a digression on this topic, see
von Bergmann \cite{vB}.) We would like to point out that both the
knot type of $\gamma$ in the fold and its framing depend on the achiral
Lefschetz fibrations we use in the construction, so does the simple model
we get. The following example illustrates this phenomenon.  \end{remark}

\begin{example} \label{FoldedSphere} \rm
If we construct $S^4$ following the recipe given in the proof of
\fullref{AllFolded}, we get $W= D^2 \times D^2 \bigcup_{C \times D^2}
D^2 \times D^2 = S^2 \times D^2$, which can be identified with $M_0$, and
the simple folded symplectic form on $Y=S^2 \times S^2=M_0\bigcup_{S^2
\times \partial D^2} M_0$ can be constructed so that its restriction
to each copy of $M_0$ is indeed the standard form $\omega_0$. Note
that here both open books already agree, so we do not need to alter
the contact structures and change the initial fibrations. Now if we
undo the surgery, that is if we trade $N=N_0$ and $M$ in the proof,
what we get is the standard folded symplectic form $\omega_0$ on $S^4$.

\labellist\small
\pinlabel {$0$} [bl] at 161 656
\pinlabel {$0$} [bl] at 243 656
\pinlabel {$0$} [l] at 258 626
\pinlabel {$\cup$ a 4--handle} [t] at 227 615
\pinlabel {$S^2{\times}S^2$} [b] at 227 580
\pinlabel {$-1$} [t] at 336 633
\pinlabel {$-1$} [l] at 354 614
\pinlabel {$0$} [bl] at 397 662
\pinlabel {$0$} [tr] at 439 605
\pinlabel {\large$\approx$} at 422 627
\endlabellist

\begin{figure}[ht!]
\centerline{\includegraphics[width=12cm]{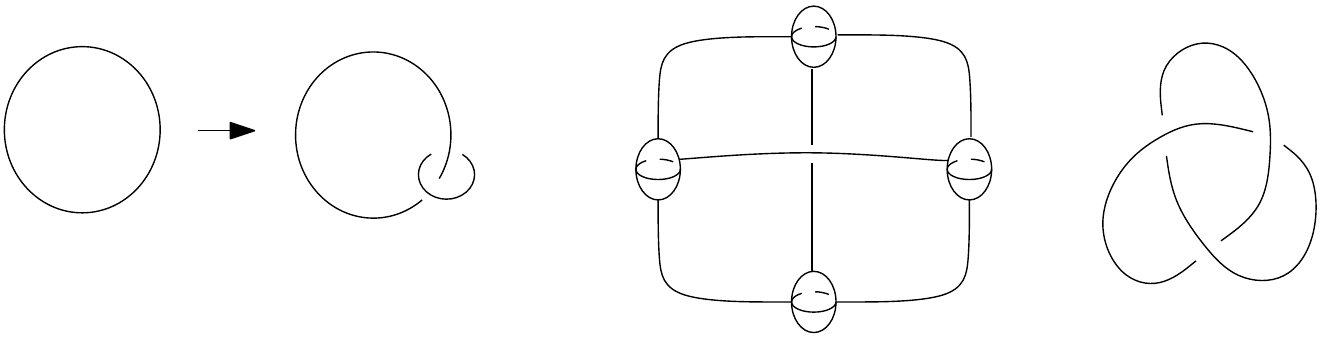}}
\caption{On the left: $0$--surgery along the binding yields a
trivial $S^2$ fibration over $D^2$ on each piece, which make up $S^2
\times S^2$. On the right: $0$--surgery along the \textit{new} binding
yields a cusp neighborhood on both sides.}
\label{FoldedS4Bindings}
\end{figure}

\labellist\small
\pinlabel {$-1$} [b] at 134 403
\pinlabel {$+1$} [t] at 134 397
\pinlabel {$-1$} [r] at 172 379
\pinlabel {$+1$} [l] at 177 379
\pinlabel {$0$} [bl] at 234 439
\pinlabel {$0$} [tl] at 257 352
\pinlabel {$\cup$ two 3--handles and a 4--handle} at 200 330
\pinlabel {\large$\approx$} at 280 400
\pinlabel {$0$} [r] at 370 420
\pinlabel {$-1$} [b] at 334 403
\pinlabel {$+1$} [t] at 334 397
\pinlabel {$+1$} [l] at 376 379
\pinlabel {$0$} [bl] at 436 439
\pinlabel {$0$} [tl] at 458 352
\pinlabel {$\cup$ two 3--handles and a 4--handle} at 400 330
\pinlabel {\large$\approx$} at 85 241
\pinlabel {$+1$} [t] at 134 242
\pinlabel {$0$} [bl] at 175 250
\pinlabel {$0$} [bl] at 236 283
\pinlabel {$0$} [tl] at 257 195
\pinlabel {$\cup$ two 3--handles and a 4--handle} at 200 165
\pinlabel {$0$} [t] at 290 230
\pinlabel {$+1$} [tl] at 320 230
\pinlabel {\large$\approx$} at 355 240
\pinlabel {$0$} [t] at 402 230
\pinlabel {$+1$} [tl] at 432 230
\pinlabel {$\cup$ a 4--handle} at 415 205
\pinlabel {$\CP^2\#\overline{\CP}^2$} at 415 175
\endlabellist
\begin{figure}[ht]
\centerline{\includegraphics[width=12cm]{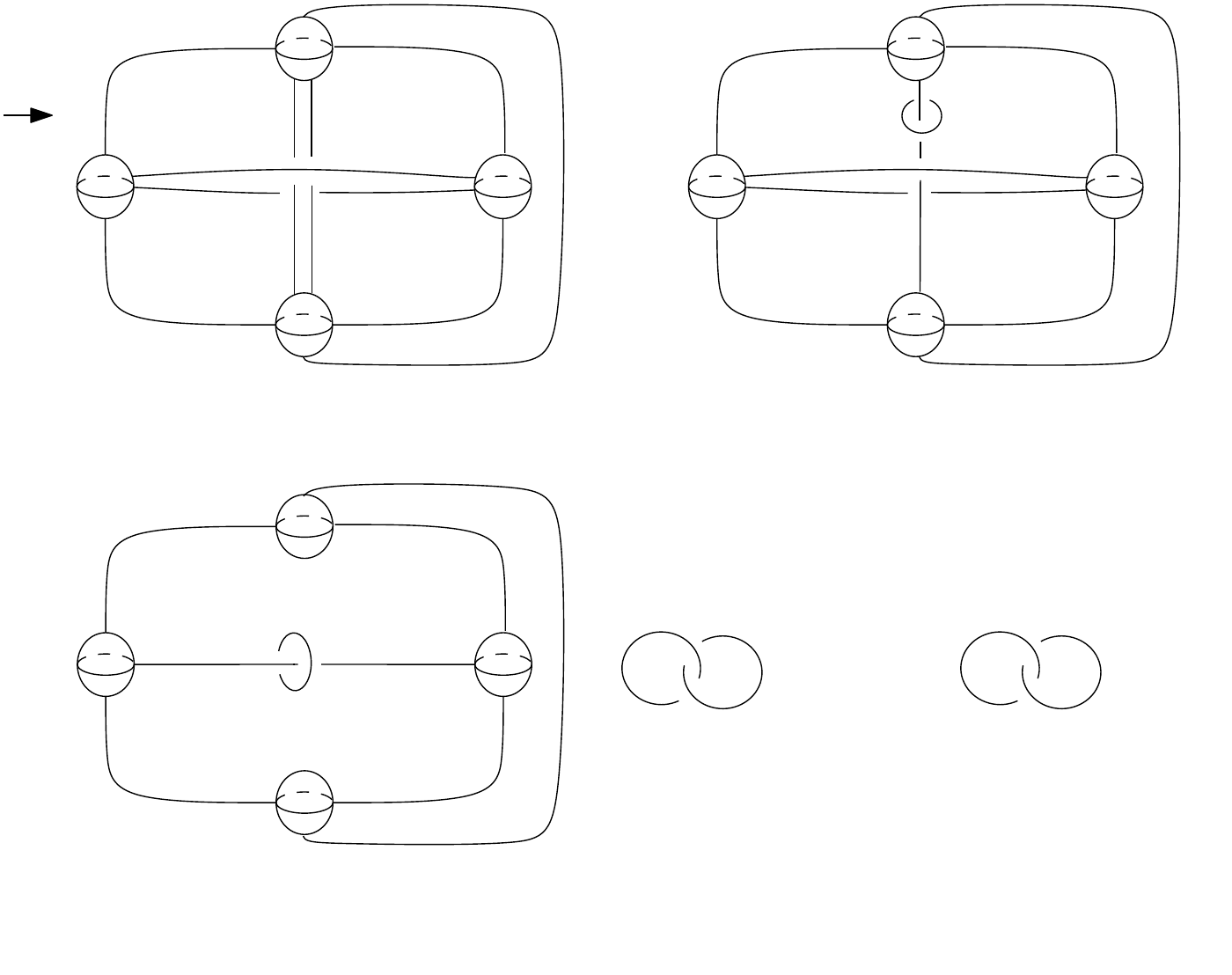}}
\caption{The achiral Lefschetz fibration on the second associated
model. The total space is shown to be $\CP^2 \# \overline{\CP}^2$.}
\label{FoldedS4SecondModel}
\end{figure}
\end{example}

It is a standard fact that surgery along a framed curve in a
simply-connected $4$--manifold will result in connect summing with either
$S^2\times S^2$ or $S^2\widetilde{\times}S^2$, depending on the framing,
which can be thought as an element of $\pi_1(\mathrm{SO}(3))=\Z_2$. In
Etnyre--Fuller \cite{EF} (also see Harer \cite{Ha}) it is described
how one can homotope the framed knot in the $3$--manifold to another
framed knot, which is isotopic in the ambient $4$--manifold to the
original one, so that their framings differ by one and that surgering
the new curve yields an achiral Lefschetz fibration on the resulting
manifold as well. Applying this trick to our example, we can instead
pass to an achiral Lefschetz fibration on $S^2\widetilde{\times}S^2
\cong \CP^2 \# \overline{\CP}^2$, which is a torus fibration with
two cusp fibers of opposite signs (Figures \ref{FoldedS4Bindings} and
\ref{FoldedS4SecondModel}). The monodromy of this achiral fibration
is $t_a \, t_b \, t_b^{-1} \, t_a ^{-1}$, and the corresponding Kirby
diagram is depicted in \fullref{FoldedS4SecondModel} (see Gompf
\cite{Go1}). To verify that this manifold is $\CP^2 \# \overline{\CP}^2$,
we first slide one of the vertical $2$--handles over the other one,
and then separate this pair from the rest of the diagram by sliding over
the $0$ framed $2$--handle. Now the rest of the diagram can be shown to
be $S^4$ after obvious handle cancelations. It is possible to see that
our construction method will give a folded symplectic structure on $S^4$
equivalent to the standard one again. (Take the first example minus the
neighborhood of the binding, do a pair of mirror stabilizations on each
side, and then proceed as in the proof of \fullref{thesame}.) Note that
the first simple model above is obtained by surgering the unknot, whereas
the second comes from surgering the right trefoil in $S^3$. Surgery
framings on them differ by one in $S^4$.


\section{Decomposition theorem} \label{Decomposition}

While we shift our attention to Stein structures, we would like to have
only non-separating vanishing cycles in our constructions, as it is
suggested by the correspondence between PALFs and compact Stein surfaces
established by Loi and Piergallini \cite{LP} and by Akbulut and Ozbagci
\cite{AO}. We start with the following lemma:

\begin{lemma}
Let $X$ be a closed oriented $4$--manifold. Then it can be decomposed
into two handlebodies, each of which admits an {\em allowable} achiral
Lefschetz fibration over $D^2$, such that the fibers have connected
boundaries and that the induced open books are the same.
\end{lemma}

\begin{proof}
We follow the construction of Etnyre and Fuller with more care given to
having fibrations allowable. Take a handlebody decomposition of $X$ with
one $0$-- and one $4$--handle, let $X_1$ be the union of the $0$--handle,
$1$--handles and $2$--handles, and $X_2$ be the union of $3$--handles
and the $4$--handle. As it was implicitly present in the thesis of Harer
\cite{Ha}, and was also observed by Akbulut and Ozbagci \cite{AO}, one can
always build an achiral Lefschetz fibration on a given $2$--handlebody
so that all vanishing cycles are non-separating. Therefore, we can
start with allowable fibrations and then proceed with stabilizations
as described by Etnyre and Fuller \cite{EF} to match the induced open
books. A stabilization is given by gluing a positive or a negative
Lefschetz $2$--handle along a new $1$--handle added to a regular fiber,
and in order to keep the binding connected, we always introduce another
adjacent stabilization. Therefore, all vanishing cycles introduced during
stabilizations are also nonseparating. The proof is completed by
induction.
\end{proof}

\begin{theorem} \label{SteinDecom}
Let $X$ be a closed oriented smooth $4$--manifold. Then $X$ admits a
decomposition into two codimension zero submanifolds $X_+$ and $X_-$,
such that $X_+$ and $-X_-$ are both compact Stein manifolds with strictly
pseudoconvex boundaries. These Stein structures can be chosen so that
the induced contact structures $\xi_+$ on $\partial X_+$ and $\xi_-$ on
$-\partial X_-$ are isotopic. Furthermore, there are PALFs on each piece
such that the open book decompositions they induce on $\partial X_+$ and
$-\partial X_-$ are compatible with $\xi_+$ and $\xi_-$, respectively,
and they coincide. In short, all data match on the hypersurface $H=
\partial X_+ = -\partial X_-$.
\end{theorem}

\begin{proof}
The lemma above gives us a decomposition of $X$ into two pieces $X_1$ and
$X_2$ with allowable achiral Lefschetz fibrations $f_1$ and $f_2$ on them,
such that induced open books on the boundaries match. As in the proof
of \fullref{AllFolded}, we glue these two pieces along the truncated
pages to get
$$W= X_1 \bigcup_{f_1^{-1}(\partial D^2)=f_2^{-1}(\partial D^2)} X_2.$$
Next we glue in $S^2 \times D^2$ to cap off the fibers
and establish an achiral Lefschetz fibration $\hat{f}\co Y \to S^2$
with closed fibers.

We wish to split the base of this fibration into two disks $D_+$ and $D_-$
so that all the positive critical values lie in the interior of $D_+$
and all the negative ones lie in the interior of $D_-$. As discussed
earlier, restrictions of $f$ give a positive Lefschetz fibration on $X_+
= f^{-1}(D_+)$ and a negative Lefschetz fibration on $X_- = f^{-1}(D_-)$,
respectively. It also prescribes a surface bundle over $S^1 = \partial D_+
= -\partial D_-$ on the hypersurface separating $Y_+$ and $Y_-$. This
time we would like to describe the splitting more carefully by taking
into consideration how the global monodromies of the new fibrations are
related to the original ones.

\labellist\tiny
\hair=1pt
\pinlabel {$q_3^1$} [tr] at 178 384
\pinlabel {$q_2^1$} [r] at 200 374
\pinlabel {$q_1^1$} [r] at 218 353
\pinlabel {$q_1^2$} [tr] at 219 323
\pinlabel {$q_2^2$} [br] at 209 290
\pinlabel {$q_3^2$} [r] at 168 269
\pinlabel {$q_{k_2}^1$} [tl] at 124 308
\pinlabel {$q_{k_1}^1$} [tl] at 122 341
\pinlabel {\large$D$} [tl] at 246 275
\pinlabel {$q_3^1$} [tr] at 400 384
\pinlabel {$q_2^1$} [r] at 422 374
\pinlabel {$q_1^1$} [r] at 440 353
\pinlabel {$q_1^2$} [tr] at 451 323
\pinlabel {$q_2^2$} [br] at 430 290
\pinlabel {$q_3^2$} [r] at 390 269
\pinlabel {$q_{k_2}^1$} [tl] at 346 308
\pinlabel {$q_{k_1}^1$} [tl] at 344 341
\pinlabel {\large$D$} [tl] at 468 275
\endlabellist
\begin{figure}[ht]
\centerline{\includegraphics[width=12cm]{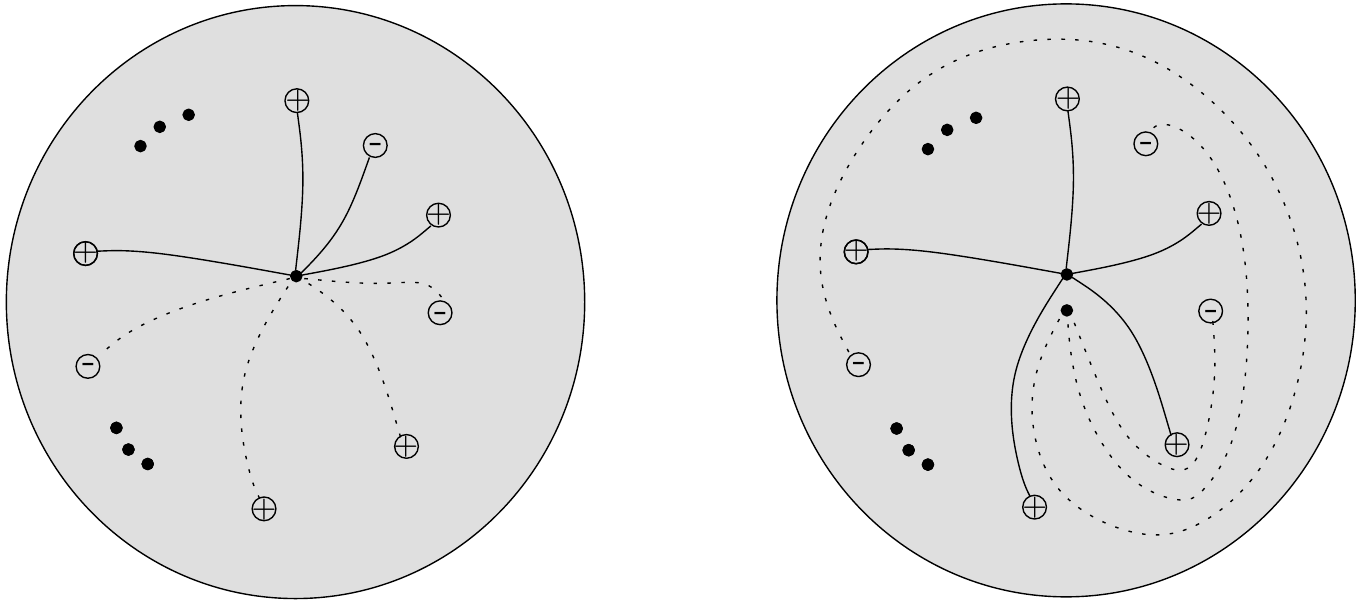}}
\label{Monodromies}
\caption{New monodromies from old ones. On the left: $\mu_1$ is
given by solid arcs, and $\mu_2$ by dotted ones. On the right: Solid arcs
are the \textit{positive arcs} representing $\mu_+$, whereas conjugated
dotted arcs are \textit{negative}, providing a representation of $\mu_-$
after closing the base back to $S^2$.}
\end{figure}

Let $\mu_1$ be the monodromy of the achiral Lefschetz fibration $f_1$
on $X_1$ and $\mu_2$ be the monodromy of the fibration $f_2$ on
$X_2$. Fix a representation of $\mu_1$ by using arcs $\smash{s^1_1, \ldots,
s^1_{k_1}}$ and a representation of $\mu_2$ by using arcs $\smash{s^2_1, \ldots,
s^2_{k_2}}$. Corresponding critical values are denoted by $\smash{y^1_i}$
and $\smash{y^2_j}$. Monodromies of the open book decompositions bounding
each fibration are given by $\mu_1$ and $\mu_2$ as well, and they
coincide under an orientation reversing diffeomorphism, so $\mu_1 =
(\mu_2)^{-1}$. Let $V$ be a small neighborhood of a regular value
in the base $S^2$. We obtain an achiral Lefschetz fibration
$$f\co W \setminus f^{-1}(V) \longrightarrow D = S^2 \setminus V,$$
which closes up to a fibration over $S^2$. If $g$ is the genus of the page,
then this fibration is determined by the relation $\mu_1 \mu_2 = 1$ in
$\Gamma_{g,1}$ and is mapped under the maps $\Gamma_{g,1} \to \Gamma_g^1
\to \Gamma_g$ to the relation that describes the achiral Lefschetz
fibration
$$\hat{f}\co Y \setminus \hat{f}^{-1}(V) \to D.$$
Since this map
factors through $\Gamma_g^1$, the achiral Lefschetz fibration $\hat{f}$
comes naturally with a section $S$ of self-intersection zero. We denote
images of the elements in $\Gamma_{g,1}$ under this map with the same
elements, so $\mu_1 \mu_2 =1$ is the global monodromy of the achiral
Lefschetz fibration on $Y\setminus \hat{f}^{-1}(V)$. Note that we can
use the same arcs $s^1_i$ and $s^2_j$ to represent the global monodromy
of this fibration. Now, if we move counterclockwise and choose only the
arcs that run through positive Dehn twists, we establish a monodromy
$\mu_+$. Call these arcs \emph{positive}. Next, we choose a nearby base
point, and move counterclockwise by running through the negative Dehn
twists only, while avoiding intersecting any positive arc. This way,
we obtain a monodromy $\mu_-$. The new set of arcs involved in this
monodromy will be referred as \emph{negative} arcs. Observe that each
negative arc is obtained by traveling around some old arcs $s^1_i$
and $s^2_j$ in order to avoid intersecting positive arcs, then going
around the aimed negative critical point once, and finally going all
the way back on the same detour (Figure 3). That is, each negative arc
corresponds to a conjugate of a negative Dehn twist in $\Gamma_g$, which
defines a negative Dehn twist, too. By taking regular neighborhoods of
these arcs such that positive and negative arcs stay apart, we get a
disk enclosing only positive critical points, and an annulus containing
only negative critical points. Closing the fibration to a fibration over
$S^2$, the latter becomes a disk as well. Now we can enlarge any one
of these disks so both disks share the same boundary, and call the one
containing positive values $D_+$, and the other one $D_-$. So we have
a new factorization of the global monodromy of $\hat{f}$, given by the
relation $\mu_+ \mu_- =1$. The section $S$ prescribes how to lift the
new elements $\mu_{\pm}$ of $\Gamma_g$ to $\Gamma_g^1$ uniquely.

We proceed with taking out the tubular neighborhood $S^2 \times D^2$ of
the section from $Y = Y_+ \cup Y_-$, and we get an inherited splitting
$X_+ \cup X_-$. The discussion above shows that $\mu_+$ defines a positive
Lefschetz fibration on $X_+$ and $\mu_-$ defines a negative fibration
on $X_-$. To recover the original $4$--manifold $X$ we need to put back
in $S^1 \times D^3$, which has the same effect as gluing each other
the tubular neighborhoods $S^1 \times D_+^2$ and $S^1 \times D_-^2$
of the bindings of open books on $\partial X_+$ and $\partial X_-$,
respectively. Therefore we can think of $X$ as decomposing into $X_+$ and
$X_-$. We claim that this decomposition possesses the desired properties.

When we take out a tubular neighborhood of $S$ from $Y$, we turn the
positive and negative Lefschetz fibrations on $Y_+$ and $Y_-$ into a
PALF on $X_+$ and a NALF on $X_-$, respectively. In the meantime the
surgery converts the surface fibration that separates $Y_+$ and $Y_-$
to an open book decomposition on the common boundary $H=\partial
X_+=-\partial X_-$. The binding of this open book is the identified
bindings of $\partial X_+$ and $-\partial X_-$, the page $F$ is the
bounded closed surface obtained by cutting off a disk from the regular
fiber of $\hat{f}$, and the monodromy is induced from the fibration on
either side. Noting that the NALF on $X_-$ becomes a PALF on $-X_-$,
we see that both PALFs induce the same open book decomposition on their
boundaries.

By \fullref{PALF}, both $X_+$ and $-X_-$ admit Stein structures. We will
construct these Stein structures using Eliashberg's characterization so
that they match on the common boundary. The technique we are going to
use is the same as the one which was presented by Akbulut and Ozbagci
\cite{AO}: The PALF on $X_+$ is obtained by attaching positive Lefschetz
handles $h_1, \ldots, h_m$ to $X_0= F \times D^2$, which has the obvious
PALF defined by projection onto $D^2$ component. The same is true
for the PALF on $-X_-$. $F \times D^2$ has a natural Stein structure
by \fullref{Eliashberg}. We can assume all vanishing cycles (coming
from either side) sit in various pages of the open book on $H$. Read
backwards, we can think of the fibrations as being constructed by
attaching positive and negative Lefschetz handles to $H$ on either
side in a sequence following the monodromy of the open book. Thus we
can induct on the number of handles. Assume that the PALF on $X_{i-1} =
X_0 \cup h_1 \cup \ldots \cup h_{i-1}$ \, ($i\leq m$) induces an open book
decomposition on its boundary, and it carries a Stein structure such that
the contact structure induced on the boundary is compatible with this open
book. Let $C$ be the vanishing cycle of the positive Lefschetz handle
$h_i$ lying on a page $F$ of $\partial X_{i-1}$. We open up the open
book decomposition and choose a page against $F$, and glue them together
along the binding $B$ to form a smooth closed convex surface $\Sigma$
in the $3$--manifold $\partial X_{i-1}$. As $C$ is non-separating,
$\Sigma \setminus C$ is connected and it contains the dividing set,
namely $B$. So we can use the \dfn{Legendrian realization principle}
(see Giroux \cite{Gi1} and Honda \cite{Ho}) to isotope $\Sigma$ through
convex surfaces to make $C$ \, Legendrian. Note that this is done by
a small $C^{\infty}$ isotopy of the contact structure supported in a
neighborhood of $\Sigma$, which fixes the binding pointwise. Hence the
framing of $C$ relative to the fiber $F$ is the same as its contact
framing, implying that the Lefschetz handle $h_i$ is attached along
a Legendrian curve with framing $tb - 1$. By \fullref{Eliashberg}
the Stein structure extends over this handle, and as shown by Gay in
\cite{Ga}, the new open book on $\partial X_i$ will be compatible with
the new induced contact structure on $\partial X_i$. This completes
the induction. Repeating the same argument dually for $-X_-$, we see
that the compatible open books on $\partial X_+$  and $\partial(-X_-)$
are isotopic, and therefore the induced contact structures $\xi_+$
on $\partial X_+$ and $\xi_-$ on $\partial(-X_-)=-\partial X_-$ are
isotopic as well. So we fulfill all the matching conditions listed in
the statement of the theorem.
\end{proof}

\begin{remark}
Akbulut and Matveyev \cite{AM} asked if one could decompose a given
closed oriented smooth $4$--manifold into Stein pieces so that the induced
contact structures on the separating $3$--manifold coincide. Our theorem
gives an affirmative answer to this question. In the same article authors
remark that it is possible to alter their Stein decomposition to make
the induced contact plane distributions homotopic, but the tightness
of the contact structure precludes the use of Eliashberg's celebrated
theorem on overtwisted contact structures to conclude more. Considering
the underlying PALFs and isotopies of open books gives a way around this
difficulty, thanks to Giroux's Theorem.
\end{remark}

\begin{remark} \rm \label{cobordism}
In \cite{Q}, Quinn studied so-called \dfn{dual decompositions} of
$4$--manifolds: descriptions of $4$--manifolds as a union
of two $2$--handlebodies. The author formulates the same question as
in \cite{AM} in terms of necessary sequence of Kirby moves to relate
a possibly nonmatching Stein decomposition. \fullref{SteinDecom}
provides an \textit{implicit} answer to this question, and we would like
to take this as an opportunity to summarize the handle calculus behind
our construction: An arbitrary Stein decomposition $X= X_1 \cup X_2$ comes
with some PALF pair. Using the stabilization moves of Etnyre and Fuller,
we first change this PALF pair with a matching pair. This corresponds to
adding canceling $1$-- and $2$--handles to each $X_i$, or in
other words, we add canceling handle pairs of index 1, 2 and 3
in the original handlebody decomposition of $X$. In the next step, we
pass to a cobordant $4$--manifold $Y$ so that we can split the positive
and negative Lefschetz handles. Then we `undo' the surgery and get the
decomposition $X= X_+ \cup X_-$ with Stein structures on each piece that
coincide on the common boundary. Having the simply-connected case in
mind, this intermediate step can be seen as a stabilization. Let $W\cong
X\setminus S^1 \times D^3$ be the complement of a regular neighborhood
of the framed knot $\gamma$ in $X$, then the first surgery defines
a cobordism
\[
[0,2] \times W \bigcup_{[0,2] \times S^1\times S^2} \bigl([0,1] \times
S^1 \times D^3 \bigr) \cup_{1 \times S^1 \times S^2} \bigl([1, 2] \times D^2
\times S^2\bigr),
\]
which is identity on the first component. We trade $2$--handles of
$X_1$ and $X_2$ in $Y$ by making use of the two extra handles of index
$2$. Finally, the composition of two cobordisms that gives back $X$ can
be seen as the double of the cobordism above, and thus it deformation
retracts to
\[
W \bigcup_{S^1\times S^2} \bigl(S^4\cup_{D^2 \times S^2} S^4\bigr).
\]
This cobordism is built by attaching cells to $\partial W = S^1 \times
S^2$, where $D^2 \times S^2$ is attached uniquely and the framing of
$\gamma$ indicates in which one of the two ways we shall glue the other
two $S^1 \times D^3$ pieces. Although here we started with a (nonmatching)
Stein decomposition, it is clear that the same discussion can be carried
out in our main construction as well. Therefore we have a well-defined
process, during which we first inflate the number of handles in a given
decomposition of $X$, and then trade some of the $2$--handles through a
cobordism to achieve the desired decomposition at the end.
\end{remark}


\section{Folded K\"ahler structures} \label{FoldedKahler}

Unlike symplectic structures, random folded symplectic structures do
not need to bear any information about the geometry or topology of the
manifold they are defined on. In order to specify more meaningful members
of this family, one first of all needs to impose some boundary conditions
on the folding hypersurface. We would like to acknowledge a result of
Kronheimer and Mrowka: In \cite{KM}, the authors prove that a compact
symplectic $4$--manifold $(Y, \omega)$ with strictly pseudoconvex boundary
has $SW_Y( K ) = 1$, where $K$ is the canonical class of $\omega$. This
motivates us to see such manifolds as building blocks of $4$--manifolds,
and yields a good boundary constraint for folded symplectic structures,
at least in this dimension. Henceforth, we assume that the fold
$H=(\omega^n)^{-1}(0)$ of a given folded symplectic manifold $(X^{2n},
\, \omega)$ is always connected and nonempty. We will generalize the
notion of a K\"ahler structure on a smooth $2n$--manifold by considering
a distinguished subset of the family of folded symplectic structures,
and we then present some properties of these structures:

\begin{definition} \label{FoldedConvex}
A folded symplectic form $\omega$ on an oriented $2n$--dimensional manifold
$X$ is called a \dfn{folded K\"ahler structure}, if there is a tubular
neighborhood $N$ of $H$ such that:
\begin{enumerate}
\item The closure of each component of $X \setminus N$ is a compact
K\"ahler manifold
$$\bigl(\pm X_{\pm}, \omega|_{\pm X_{\pm}}\bigr)$$
with strictly pseudoconvex boundary,
\item $(N, \omega|_N)$ is folded symplectomorphic to $(\,[-1, 1] \times
H , \, d((t^2+1) \, \pi^*(\alpha) \,)$, where $\alpha$ is a contact
$1$--form on the fold $H$, and $\pi$ is the projection
$\pi\co[-1,1]{\times}H{\to}H$.
\end{enumerate}
In addition, if each $(\pm X_{\pm}, \omega|_{\pm X_{\pm}})$ is strictly
pseudoconvex, we say $\omega$ is a \dfn{nicely folded K\"ahler structure}
on $X$.
\end{definition}

In the definition above, \dfn{nice folding} can be reformulated as folding
Stein manifolds along matching strictly pseudoconvex boundaries. Recall
that if $\psi$ is a proper strictly plurisubharmonic function on a complex
manifold $S$, then the associated $2$--form $\omega_{\psi} = - dJ^*d\psi$
is K\"ahler, and importantly, the symplectic class of $(X, \omega_{\psi})$
is independent of the choice of $\psi$ (see Eliashberg and Gromov
\cite{ElGr}). Therefore, to
complete our alternative formulation, we ask that each piece $\pm
X_{\pm}$ should admit some proper strictly plurisubharmonic function
$\psi_{\pm}$, so that $\omega|_{\pm X_{\pm}}=\omega_{\psi_{\pm}}$. In
short, it is built in the definition that a nicely folded K\"ahler
manifold is folded K\"ahler. Finally note that, due to a theorem of
Bogomolov \cite{B}, any compact folded K\"ahler manifold $X$ can be made
nicely folded after deforming the complex structure and blowing down any
exceptional curves. Even though these definitions narrow the family of
folded symplectic structures quite a lot, it is important to note, at
least in dimension four, that we still have an adequately large family
in the light of the following result:

\begin{theorem} \label{AllKahler}
Any closed oriented $4$--manifold $X$ admits a nicely folded K\"ahler
structure.
\end{theorem}

\begin{proof}
By \fullref{SteinDecom}, $X$ can be decomposed into two compact Stein
manifolds $X_+$ and $-X_-$ with strictly pseudoconvex boundaries such
that both induce the same contact structure on the common boundary $H=
\partial X_+=-\partial X_-$. We begin with adding collars $\pm U_{\pm}$
to $(\pm X_{\pm}, \omega_{\pm})$, and extending the symplectic structures
to $\omega'_{\pm}$ on $\pm X'_{\pm} = \pm (X_{\pm} \cup U_{\pm})$ so
that new boundaries $\partial (\pm X'_{\pm})$ are still convex and
contactomorphic. Let $\xi_{\pm}$ be the induced contact structures
on $\partial (\pm X'_{\pm})$ and $\psi$ be a contactomorphism between
them. Using the symplectic cut-and-paste argument of Etnyre \cite{E},
we can add a symplectic collar to $(\partial X'_+, \omega'_+)$ so that
the new boundary is not only contactomorphic to $-\partial X'_-$ but
also the induced \textit{contact forms} agree up to a multiple $k \in
\R^+$. For the sake of brevity, let us assume that $U_+$ above contains
this collar part as well. So after rescaling $\omega'_-$ (and $\omega_-$)
by $k$ if necessary, we see that restrictions of symplectic forms
$\omega'_+|_{\partial X'_+}$ and $k \omega'_- |_{-\partial  X'_-}$ agree
via $\psi$, and orientations of null-foliations (which correspond to Reeb
directions) are the same. Therefore, once again we can apply the folding
technique of Cannas da Silva, Guillemin and Woodward \cite{CGW} to obtain a folded symplectic structure $\omega$
on $X'_+ \cup X'_-$ such that $\omega$ agrees with $\omega'_{\pm}$ on the
complement of a small tubular neighborhood of the fold $H$. We enlarge
this neighborhood to include $U_+$ and $U_-$ and call it $N$. It follows
that $X= X_+ \cup N \cup X_- \cong X_+ \cup X_-$, and $\omega |_{X_+}
= \omega_+$, whereas $\omega |_{X_-} = k \omega_-$. Also note that,
the folding operation provides us with the desired local model on $N$,
that is, $(N, \omega|_N)$ is folded symplectomorphic to $(\,[-1, +1]
\times H , \, d((t^2+1) \, \pi^*(\alpha) \,)$ by construction \cite{CGW}.

Lastly, suppose $\psi_{\pm}\co \pm X_{\pm} \to \ [0, \infty)$ are proper
strictly plurisbuharmonic functions such that $\pm \partial X_{\pm}$
correspond to the maximum points of $\psi_{\pm}$, and $\omega_{\pm} =
-dJ^*d\psi_{\pm}$, respectively. If $k \neq 1$, we can replace $\psi_-$
with $k \psi_-$ and obtain $k \omega_-$ above as a K\"ahler form of a
strictly pseudoconvex manifold. Equipped with these properties, $\omega$
is a nicely folded K\"ahler form on $X$.
\end{proof}

\begin{remark} \rm \label{thesame}
It is clear that \fullref{AllKahler} is a refinement of \fullref{AllFolded}. Since the folded forms we have constructed in both proofs
are obtained through similiar steps, one expects that these structures
are actually equivalent. Next, we verify this fact, and this way we
get an insight of how folded K\"ahler forms are `supported' (precise
definition is given below) by Lefschetz fibrations as was illustrated
in \fullref{SimpleModel}:

Take the PALF on $X_+$ in the proof of \fullref{AllFolded}, and attach
a symplectic $2$--handle along the binding of the induced open book on
$\partial X_+$ as described by Eliashberg in \cite{El2}. This yields a
symplectic Lefschetz fibration over $D_+^2$. Dually the same argument
for the NALF on $X_-$ gives an anti-symplectic Lefschetz fibration
over $D_-^2$, and these handle attachments can be done so that two
fibrations agree on the common boundary. Moreover, we can assume that
these fibrations have genus at least two, so the fibrations can be matched
as symplectic surface fibrations over a circle, as it was pointed out
in \cite{El2}. At the end we get a \textit{simple folded symplectic}
manifold $Y$ obtained from $X$ after a surgery along a framed curve
$\gamma$. However, any two simple folded symplectic forms compatible
with a fixed achiral Lefschetz fibration are deformation equivalent by
\fullref{SimpleModel}. Moreover, we can normalize both forms
on the disks which are parallel copies of cocores of new $2$--handles
that were used to cap off the fibers. Hence, these two folded forms are
deformation equivalent on $Y \setminus S^2 \times D^2$. As the folded
symplectic structure on $D^3 \times S^1$ which is glued back in to
recover $X$ is standard, the folded symplectic form constructed in the
proof of \fullref{AllFolded} and the folded K\"ahler form obtained in
\fullref{AllKahler} are indeed \emph{equivalent as folded symplectic
structures}.
\end{remark}

Motivated by symplectic and near-symplectic cases (see Donaldson
\cite{D} and Auroux--Donaldson--Katzarkov \cite{ADK}), we can conclude
our discussion above by defining the Lefschetz fibration analogue for
our structures:

\begin{definition}
Let $X$ be a closed oriented $4$--manifold, and decompose $S^2$ as
the union of the upper-hemisphere $D_+$ and the lower-hemisphere $D_-$
which are glued along the equator $C=\partial D_+ = -\partial D_-$. Then
a smooth map $f\co X \to S^2$ is said to be a \dfn{folded Lefschetz
fibration} on $X$, if it restricts to a PALF over $D_+$, to a NALF over
$D_-$, and to an open book over $C$ bounding both fibrations.
\end{definition}

\begin{definition} \label{CompatibleFolded}
Let $X$ be a closed oriented $4$--manifold equipped with a nicely folded
K\"ahler form $\omega$. Then a folded Lefschetz fibration $f\co X \to S^2$
is said to be \dfn{compatible with $\omega$} if each Stein piece $X_{\pm}$
corresponds to $f^{-1}(D_{\pm})$, and if the contact structure they
induce on $H=f^{-1}(C)$ is compatible with the open book decomposition
coming from the fibration $f$. In this case, we also say that nicely
folded K\"ahler manifold $(X, \omega)$ is \dfn{supported by the folded
Lefschetz fibration $f$}.
\end{definition}

The compatibility in the above definition is completely \emph{on the
symplectic level}. This becomes more visible if once again we recall
that surgering the binding $\gamma$ of the open book $\smash{f|_{H \setminus
\gamma}\co H \setminus \gamma \to S^1}$, we pass to a \emph{simple
model} where the folded Lefschetz fibration can be extended to a folded
\emph{symplectic} achiral Lefschetz fibration $\hat{f}$ with closed
fibers. Also note that, since Stein manifolds harbor less topological
obstructions in complex dimensions $> 2$, it is very likely that they
admit higher dimensional analogues of PALFs with similar topological
correspondences. If that is established, last two definitions, as well as
several results in this paper, can be generalized to all $2n$--dimensions.

The complete statement of \fullref{SteinDecom} combined with \fullref{AllKahler} shows that, given a closed oriented $4$--manifold $X$,
one can always find a nicely folded K\"ahler structure $\omega$ on $X$
together with a compatible folded Lefschetz fibration. Next, we prove
that this property in fact holds for any nicely folded K\"ahler structure:

\begin{proposition} \label{CompatibleFoldedLF}
Any nicely folded K\"ahler structure $\omega$ on $X$, up to orientation
preserving diffeomorphism, admits a compatible folded Lefschetz fibration.
\end{proposition}

\begin{proof}
Each Stein piece $X_+$ and $-X_-$ admits a PALF by \fullref{PALF}. If
we construct these fibrations following the algorithm of Akbulut and
Ozbagci \cite{AO} and keep track of the associated open books, the work
of Plamenevskaya \cite{P} shows that we can establish PALFs $f_{\pm}\co
\pm X_{\pm} \to D^2$ with the property that the open book decomposition
on the boundaries are compatible with the contact structures induced
from the Stein structures on $\pm X_{\pm}$, respectively. As the contact
structures are assumed to be the same, \fullref{Giroux} tells us that we
can match these open books after positive stabilizations. Consequently,
we get a compatible folded Lefschetz fibration.
\end{proof}

\begin{remark} \rm
A folded Lefschetz fibration that supports a given folded K\"ahler
structure fails to be unique. In fact, one can find infinitely many
pairwise inequivalent such fibrations. This can be seen for example from
the construction of \cite{AO}, by using different $(p,q)$ torus knots
in their algorithm which we adopt for building our achiral Lefschetz
fibrations.
\end{remark}

\begin{example} \rm
The easiest examples are doubles. If $Y^4$ is a compact K\"ahler
manifold with strictly pseudoconvex boundary, then $X=Y \cup -Y$ is
equipped with a folded K\"ahler structure. When $Y$ is indeed Stein, we
get a nicely folded structure. The first folded structure constructed in
\fullref{FoldedSphere} is the double of standard $D^4 \subset \C^2$,
whereas the latter is a `monodromy double' of a cusp neighborhood minus
a section. Here by `monodromy double' we mean that the pieces are first
glued along the pages of the open books, and if the monodromy of the
folded Lefschetz fibration on one piece is $\mu$, then the monodromy on
the other one is $\mu^{-1}$.
\end{example}

\begin{example} \rm
There is a construction which also allows us to see the nicely
folded K\"ahler structure together with a compatible folded Lefschetz
fibration. Take a contact $3$--manifold $(H, \xi)$, and fix a positive
open book decomposition $(B, f)$ compatible with $\xi$. Different PALFs
bounding this open book describe (possibly) different Stein fillings
of $(H, \xi)$. Indeed there are examples of infinitely many pairwise
non-diffeomorphic contact $3$--manifolds each of which admit infinitely
many pairwise non-diffeomorphic Stein fillings constructed this way
(see Ozbagci and Stipsicz \cite{OS1}). Thus for every pair of PALFs
$X_1$ and $X_2$ that fill the same open book, we can construct a nicely
folded K\"ahler form on $X=X_1 \cup -X_2$, as designated in the proof
of \fullref{AllKahler}.
\end{example}

\begin{example} \rm
The main steps of our construction are depicted in the following simple,
albeit instructive example: Start with classical handlebody decomposition
of $X =\#_8 S^2 \times S^2$ with one $0$--handle, sixteen $2$--handles, and
a $4$--handle. Let $X_1$ be the union of $0$--, $2$-- handles, and $X_2$ be
the $4$--handle. Each piece admits a $D^2$ fibration over $D^2$. However we
wish to construct allowable fibrations, so we introduce two canceling 1--
and 2--handle pairs and two canceling 2-- and 3--handle pairs in the
original handlebody decomposition of $X$. We start building the fibrations
from the scratch: Add the $1$--handles to the $0$--handle and $3$--handles
to the $4$--handle. Attach the two canceling $2$--handles with
framing $-1$ to the union of the $0$-- handle and $1$-- handles. Attach the
other two the same way to the handlebody $X_2$, which is the union of $3$-
handles and the $4$-- handle. To simplify our description, we will label
the $1$--handles of the first handlebody as $a$ and $b$, which generate
$\pi_1$ of the torus fiber with one boundary component, and we do the
same for the $1$--handles of $X_2$ under the obvious identification. So
we obtain two achiral Lefschetz fibrations over disks with bounded torus
fibers; one with monodromy $t_a^{-1} \, t_b^{-1}$, and one with $t_b \,
t_a$. One can verify by Kirby calculus that each time we insert a pair of
Lefschetz handles prescribed by $t_a \, t_a^{-1}$ or $t_b \, t_b^{-1}$,
we introduce an $S^2 \times S^2$. (See \fullref{FoldedS4SecondModel},
and observe that here we slide off the $2$--handle pair over a $+1$
framed $2$--handle instead.) Doing this consecutively, we attach all the
remaining $2$--handles to the first handlebody, and obtain an achiral
Lefschetz fibration on $X_1$ with monodromy
\[
\mu_1 =t_a^{-1}t_b^{-1}t_b \, t_a \,t_b \,t_a \,t_b \,t_a \,t_b \,t_a
\,t_a^{-1}t_b^{-1}t_a^{-1}t_b^{-1}t_a^{-1}t_b^{-1}t_a^{-1}t_b^{-1},
\]
whereas $X_2$ still has the monodromy
\[
\mu_2 =t_b \, t_a = (t_a^{-1} \, t_b^{-1})^{-1}.
\]
Both open books that bound these fibrations contain negative Dehn twists
(recall that on $- \partial X_2$, the monodromy is $\mu_2^{-1}$), and
therefore the contact structures they support are overtwisted. As we have
already manipulated the monodromy that way, contact structures and open
books are isotopic, so we can glue $X_1$ and $X_2$ along the truncated
pages. Putting in $S^2 \times D^2$ we pass to a torus fibration $\hat{f}
\co Y \to S^2$ with global monodromy $\mu_1 \cdot \mu_2$. (Applying
the handle slides given in \fullref{FoldedSphere} repeatedly, and
proceeding with the same handle cancelations, one can indeed check that
$Y= \#_8 S^2 \times S^2 \# \CP^2 \# \overline{\CP}^2$.) Now the monodromy
splits easily as explained in the proof of \fullref{SteinDecom}, and
we get $\mu_+=(t_b \, t_a)^5$ and $\mu_-=(t_a^{-1} \, t_b^{-1})^5$. It is
not hard to see that when we take out the section now, we get pieces $X_+$
and $X_-$, which are diffeomorphic to $-E_8$ and $E_8$, respectively. So
$X$ decomposes into a Stein piece $-E_8$ and an anti-Stein piece
$E_8$. This defines a nicely folded K\"ahler structure $\omega$ on $X$,
folded along the Poincar\'e homology sphere $\Sigma(2,3,5)$, and
it is supported by a folded Lefschetz fibration which is the monodromy
double of a torus fibration over $D^2$ minus a section.
\end{example}


\bibliographystyle{gtart}
\bibliography{link}

\end{document}